# Mentoring Undergraduate Interdisciplinary Mathematics Research Students: Junior Faculty Experiences

*Jana L. Gevertz*
Department of Mathematics & Statistics
The College of New Jersey
2000 Pennington Road
Ewing, NJ, 08628, USA
gevertz@tcnj.edu

*Peter S. Kim*
School of Mathematics and Statistics
University of Sydney
Sydney, NSW 2006, Australia
pkim@maths.usyd.edu.au

*Joanna R. Wares*
Department of Mathematics and Computer Science
University of Richmond
28 Westhampton Way
Richmond, VA, 23173, USA
jwares@richmond.edu

**Abstract:** To be successful, junior faculty must properly manage their time in the face of expanding responsibilities. One such responsibility is supervising undergraduate research projects. Student research projects (either single or multi-student) can be undertaken as a full-time summer experience, or as a part-time academic year commitment. With many potential undergraduate research formats, and with different types of students, junior faculty may find challenges in forming their research group, establishing a structure that promotes student productivity, picking an appropriate project, or in effectively mentoring their students.

This article draws from the authors' experiences to help junior faculty navigate these complexities so that all parties reap the benefits of undergraduate research in interdisciplinary mathematical disciplines.

## 1 INTRODUCTION

For junior faculty, supervising undergraduate researchers can be truly rewarding, allowing one to teach and mentor a small group of motivated and intelligent students. Working with individuals or small groups in this way can foster close ties among those involved, as well as with the department and school at large. Students not only gain a mastery of the ideas involved in the research project, but learn to discover and reason. In addition, they learn to share their ideas, construct logical arguments, and think critically – these skills apply to mathematics, science, and life in general. In addition, faculty members can use these experiences to guide their teaching and when done exceptionally, the undergraduate students also help advance their faculty mentor's research agenda.

More mundanely, junior faculty members may improve their tenure portfolios by supervising undergraduate research projects, particularly at primarily undergraduate institutions. Some schools or particular departments of schools even require junior faculty members to supervise undergraduate student research projects before tenure. That said, this varies greatly from institution to institution and junior faculty should be aware of how undergraduate research is valued in their department.

Additionally, postdoctoral fellows and graduate students aiming for jobs as tenure-track faculty members are wise to seek student mentoring and advising opportunities. Many institutions, especially primarily undergraduate ones, value undergraduate research as an integral part of developing erudite students. These schools look to hire teacher-scholars – those who integrate their research and their pedagogy. Having experience with mentoring undergraduate research projects is one measure of being a teacher-scholar, and therefore can help one to land a tenure-track position at an undergraduate-focused institution.

While undergraduates are now active participants in many scholarly fields, interdisciplinary mathematical projects in particular lends themselves to undergraduate research because problems can be found that arise naturally from the students' own interests in the world around them. For instance, mathematical medicine, a subfield of mathematical biology, appeals to undergraduate students who frequently hear about health and disease in the media or in their personal experiences. For this reason, the students often have an external motivation to better understand and conquer the problem at hand.

In this work, we describe opportunities and challenges, based on our experiences, for new faculty members working with undergraduate students on interdisciplinary mathematical research. Our examples come from the field of mathematical biology, but much of the advice that we give will generalize to other interdisciplinary mathematical pursuits. Two of us, JLG and JRW work at liberal arts colleges (The College of New Jersey and the University of Richmond, respectively), and the third, PSK, works at a large research university (University of Sydney), where undergraduates frequently complete honors research projects. Here, we give tips based on what we have learned from our experiences supervising undergraduate research students, both as postdoctoral fellows and as junior faculty members. We hope others can learn from our experiences.

## 2 INSTITUTIONAL AND DEPARTMENTAL SUPPORT

For tenure-track faculty members, institutional and departmental support is highly recommended for undergraduate research projects to be a successful and positive experience for both the students and the faculty [11]. While there is no one correct model for supporting undergraduate research, there are commonalities among various institutions that allow for successful undergraduate research. For instance, the mission of both The College of New Jersey and the University of Richmond is based upon the teacher-scholar model, and therefore infrastructure is in place for undergraduate students and faculty to engage in research during the academic year and in the summer. At The College of New Jersey, summer research stipends are available on a competitive basis through the MUSE (Mentored Undergraduate Summer Experiences) program, and at the University of Richmond the Richmond Guarantee program ensures every student a funded research or internship experience. In both cases, students are financially supported for their summer research, and faculty members receive stipends for working with these students.

There are also support mechanisms at some institutions for academic year research, and some external grants [15] are also available for this purpose. At The College of New Jersey and the University of Richmond, students can complete research during the semester through independent research/study courses. When a student registers for such a course, they are able to get credit for their research, and at least at The College of New Jersey, the faculty member gets some credit for this as well. Depending on the institution, students can receive different numbers of credits for their research based on the number of hours they commit per week. Often times, this research experience counts towards the students major or minor requirements, giving students further incentive to conduct research during the academic year.

At institutions that do provide faculty credit for undergraduate research, it is typically the case that faculty compensation is not enough to reduce the teaching load by a full course. That said, in these instances credits from undergraduate research can create more flexibility in the professors' schedule, or can result in some financial compensation at the end of the academic year. At institutions that do not count research mentoring in the faculty load (including University of Richmond), such experiences can still be valuable for faculty and their students. However, faculty must be careful to take on a reasonable mentoring load, or to seek external support (e.g. through the Center for Undergraduate Research in Mathematics [15]) to reduce one's teaching load while engaging in research mentorship.

At the University of Sydney, and throughout Australia, academic year undergraduate research is similarly supported, but through a very different mechanism: it mostly takes place through the 4th-year honors program. In Australian universities, students receive a B.Sc. in three years, but can choose to continue to a 4th-year honors program to obtain a B.Sc.(Hons). In the mathematics honors program at the University of Sydney, 60% of the final grade is determined by coursework and 40% is determined by a 60-page thesis and presentation, which are based on research done with a supervisor. The option of an additional honors year at the end of one's undergraduate studies provides a nine-month period during which motivated undergraduates can focus on one research project.

For all junior faculty members, supervising undergraduate research is a great way to be involved with students and colleagues. However, even with institutional support, undergraduate research is a big commitment, and new faculty must be careful to protect their time. While we fully support junior faculty taking on new opportunities as they appear [10], saying yes to every interested undergraduate research student is not always in the best interest of the professor. Frequently, these projects are extra work that we voluntarily supervise because we believe the experience is invaluable to our students. However, junior faculty members have research standards to meet, and undergraduate research can slow down research productivity. Therefore, one needs to be very aware of the standards to which they are being held and how undergraduate research is viewed when going up for tenure.

## 3 FORMING AND SUSTAINING AN UNDERGRADUATE RESEARCH GROUP

Depending on the size of your institution and your research area, picking students for your research group may be a question of identifying any interested students, or choosing from a larger pool of many interested qualified students. The latter problem is becoming more common, as undergraduate research has practically become a prerequisite for admission into a mathematics doctoral program. The number of students interested in interdisciplinary research pursuits is also growing, as more institutions teach the basic STEM sequence using a cross-disciplinary approach (like the Integrated Quantitative Studies program and the Science Math and Research Training program at the University of Richmond [9]). Table 1 quantifies the authors' experiences working with undergraduate research students during their careers.

| Author | Number of Undergraduates Supervised | Years in Faculty Position |
|---|---|---|
| Gevertz | 10 | 7 |
| Wares | 7 (1 supervised during postdoc) | 5 |
| Kim | 15 | 5 |

Table 1: Quantifying authors' undergraduate research mentoring experiences.

### 3.1 JLG's Experience at a Primarily Undergraduate Institution

I have largely formed and sustained my interdisciplinary research group (in mathematical biology) by attracting students from my courses. This is especially true when I teach our Mathematical Biology course, but I have had success recruiting students from other mathematics courses as well. In differential equations, I present many examples from biology and assign a final course project that requires analyzing a mathematical model of cancer treatment [4] (which happens to be my particular research emphasis). Even when teaching second semester calculus, I try to connect techniques of the course to mathematical biology. While these gestures (particularly in calculus) are small, they can pay big dividends. Of the ten students I have supervised in undergraduate research projects over the past four years, seven of those students approached me about research after taking one of my courses. In fact, two of those students came from a calculus course.

Students have also been attracted to my group because of undergraduate-accessible campus talks I have given on my research. Other students have come from academic advising. At The College of New Jersey, students are directly advised by a faculty member in their major department. If a faculty member learns that a student has interests in both mathematics and biology, they are directed to a mathematical biologist. Finally, interdisciplinary researchers are referred research students through a means not available to other mathematicians: from other departments on campus! Three of the ten students I have worked with were Biology majors (two with math minors), and a fourth student is a double major in Applied Mathematics and Biology. A fifth student was a Biophysics major, and Biomedical Engineering students have also expressed interest in my research. What is interesting about attracting this pool of students is that they often come from academic departments with a culture of undergraduate research. Freshman Biology majors, for instance, are told to start thinking about potential research mentors from the day they set foot on campus.

With a significant pool of students expressing interest in my mathematical biology research endeavors, I have been forced to learn some valuable lessons about choosing the right students for my research group, as well as deciding upon the right number of students. When managing individual student research projects during the academic year, I have found two students per semester to be ideal. Even with students working on separate projects, comparable mathematical, biological, and computer programming skills are often required. A fellow research student provides a great sounding board, and sometimes a great teacher, when a student hits a roadblock. For those students working on academic year projects, at minimum I require those students to commit to two semesters of research. In one semester, there is simply not enough time for a student to learn all the background material and make significant progress on even a simple project.

Forming and sustaining a healthy and productive research group requires picking the right students to work with. Others have posed some excellent criteria for selecting undergraduate research students [19], so we will not repeat their advice here. There are, however, some unique criteria when it comes to selecting students for interdisciplinary mathematics research, particularly for single-student projects. No matter what the student's academic major, at minimum a student must be trained in multiple disciplines. I require that my students are proficient in the application area (in my case of mathematical biology, freshman-level biology is needed), dynamical systems, and computer programming. My experience is that a student without training in one of these areas often struggles to jump-start and progress through their project, particularly since I only try to give projects that directly stem from my main research program.

That said, it is possible to engage students with a less substantial background in research, if the project is carefully designed. However, given the large number of students that tend to express interest in my research, and given the nature of projects I have had my students work on, I do not commit to working with a student without this interdisciplinary skill set until they successfully complete the necessary coursework. Even with this strict criteria, I have occasionally run into the situation of having more qualified students than I can reasonably supervise. In this case, one strategy that has worked well for me is asking potential research students to engage in a small preliminary project. While this project requires little time and effort, I am often able to assess a student's interest and ability through the response and follow-up work on this project. Another way to assess which students will fit best is to consider their excitement towards learning in your class (if you have taught the student), or to reach out to their former professors.

**3.2 JRW's Experience at a Liberal Arts College**

In forming a research group in which students are meant to work together on a single project, the process of choosing students is different. There is room to work with students who have shown enthusiasm and maturity of thinking, even if they have not yet harnessed their creative juices enough to work independently. To be successful in a group setting, students must be able to work well with others. In this case, both the academic accomplishments of each student and their personalities must be taken into account, with a balance of leadership and cooperation in mind. When possible, faculty members should form a team of students with a mix of experiences and attributes, both academically and personally. I have found that teams consisting of both underclassmen and upperclassmen work best. The underclassmen are more likely to continue the work after the summer and the upperclassmen have more experience at the start which they can share with the others. In addition, I seek students that will be strong leaders and form a group consisting of them in combination with students who are good team players.

My most recent experience supervising undergraduate students helped me to develop some best practices for forming a group of undergraduate researchers. I was fortunate enough to collaborate with a more senior colleague of mine (Dr. Kathy Hoke) on a project to develop statistical models to help diagnose mild traumatic brain injury. Students heard about our work through our classes and also through advertising within our department. With Dr. Hoke, I learned the benefit of meeting with each interested student individually to discuss their goals for the summer and the project topic. To obtain the summer fellowship from the School of Arts and Sciences, each student must write a research proposal and submit a CV. I would recommend having your students do this even if it is not required. Through our meetings and the students' writing, we determined the particular interests of each student and helped them develop ideas for their summer work. It also became obvious who would lead the group and who would be a good team player. In the end, we had a very productive summer with four very different students.

**3.3 PSK's Experience at a Large Research University**

It is not always the case that a faculty member has control over the undergraduate students in their research group. At the University of Sydney, admission to the honors program is selective, and for this reason the honors students seem to choose their supervisor more than the other way around. Fortunately, even though I could not hand-pick the students, I have worked with a steady stream of these talented and enthusiastic students who have the drive to push substantial projects forward. My experience working with some of these students will be detailed in Section 5.

# 4 ORGANIZATION OF AN UNDERGRADUATE RESEARCH PROJECT

Undergraduate students develop into erudite thinkers by bringing a research project from conception to completion. Proper and realizable expectations set early in the research collaboration (for both the faculty mentor and the students) are an indispensable part of the process.

### 4.1 Pre-Research Phase

Given that undergraduate research is a serious investment of time and energy, it is important to determine if a project is right for you and your student before committing to it. Here are tips to follow prior to committing to an undergraduate research project:

- Have an initial discussion to explore interests and project ideas. JRW likes to work with student-generated projects (as long as they are close enough to her research for her to be able to guide them), whereas JLG and PSK have found that the best projects spin directly off of their own research program. Student-generated projects have the advantage of complete student buy-in. An advantage of faculty-generated projects is that the mentor is already well-versed in the area and can answer most questions without doing much additional work.

- Have the student write a project description or proposal. Even if not mandatory for funding (as it is at the University of Richmond), having the student write the project proposal is a great starting point for any collaboration, informing both student and professor about each other's interest and knowledge.

- Many times in interdisciplinary mathematics projects, there are different mathematics topics in one project. Discuss with each student which part of the math interests them most and in what context. For group projects, try to break off a piece of the project for each individual so that students can focus on the mathematical topic of most interest to them.

- Have the students explain why they want to participate in your research. This can help you understand what expectations the student has. Further, it is important to understand your own expectations for the student in order to ensure that the research project will be productive for both of you.

**4.2 Academic Year Projects with Single Students: JLG's Perspective**

Even with a well-designed system to support undergraduate students and faculty in academic-year research, there are still logistical challenges that make academic year research different from summer research. In the summer, students are able (and in fact, are generally required) to devote their full-time efforts to research. Conversely, in the academic year, students are balancing classes, paid work, and extracurricular activities along with their research. As a result, a typical undergraduate is only available for 6 to 12 hours a week for research. On the positive side it allows for a student to sustain their research project over months/years instead of weeks [16].

Managing multiple responsibilities does create scheduling constraints that one typically does not have to worry about in the summer. Coordinating schedules between busy faculty and students to find regular overlapping time can be challenging. Coordinating multi-student projects is even more challenging, and students who may make ideal partners on a project simply may not have availability at the same time in a semester to perform their research in groups. For this reason, although there are organizations which recommend undergraduate research should always be done in groups [15], I have found it more practical to have students work individually on projects during the academic year.

Undergraduate students are accustomed to a lot of structure from their courses. Experientially, I have found it valuable for a student to identify weekly time slots for their research at the start of the semester, with at minimum half of that time being selected to overlap with their faculty mentors availability. During this time, a research space is reserved near the faculty member's office (when physically possible) to promote student-faculty interactions during the research process. Besides trying to coordinate schedules with my undergraduate collaborators to promote informal interactions, I also require one regularly-scheduled meeting each week with each of my students, and occasionally group meetings are scheduled.

In order for a student to be productive during academic year research, it is helpful that they establish some background knowledge on their project prior to the start of the semester. For instance, for a mathematical biology project, this often means having the student read biology and mathematical modeling papers that form the foundation of our project. It is important to do this prior to the semester, as for any project, especially an interdisciplinary one, it is easy to bog a student down in background reading [2]. Then, at the beginning of the semester, it is useful to take any student questions, and to have the student prepare a presentation on their background reading. The act of presenting what they read helps us identify any misconceptions and solidify their understanding.

Almost universally, an undergraduate mathematical biology researcher will be building or analyzing a model that has its foundation in a previously-studied mathematical model. Therefore, I next have the student replicate the results of their project's foundational paper. Often times, this requires that the student writes computer code and studies some mathematics they have not seen in a course. For a mathematical biology project, this provides the students with a gentle introduction to research and helps build their confidence as they work through a challenging published paper.

As the student progresses from studying background material to performing research, it is important that the student regularly documents their progress in an electronic or hand-written research notebook. I instruct my students to update their notebook weekly with goals for the week, how they attempted to reach those goals, the results of any successes, as well as any failures. As students only have a limited amount of time for research during the academic year, this makes it easier for them to pick up the project after a several day hiatus.

Finally, the nature of mathematical biology research often requires us to wait for simulation results, which can take from seconds to days to run. Students must learn how to efficiently use their time as they wait for computer-generated results. For instance, if their code only takes a few minutes to run, they should be using that time to update their research notebook. On the other hand, if the code takes hours to run, that is a great time to read some relevant papers in the biology or mathematical biology literature, or to study up on a mathematical skill relevant to the project. Multi-tasking in this way, especially during the academic year, is essential for obtaining some measurable outcomes from the research project.

### 4.3 Multi-Student/Multi-Faculty Projects During the Summer: JRW's Perspective

Having more than one undergraduate student working on a single research project during the summer has many benefits. The most obvious benefit is that the students will have a peer with whom they can share ideas. Furthermore, together they can experience the joys and difficulties of learning how to successfully develop a research project. Strong students working in groups can develop their leadership skills and can learn how to organize other people to accomplish the task at hand. Weaker students who are not ready for independently bringing a research project from foundation to finish can contribute to a group, gain experience, and strengthen their problem solving and research skill sets.

Having more than one professor work on a single project can also benefit both the students and the professors. As a junior faculty member, working with a more senior faculty member can help you learn how to best structure your undergraduate research program and how to properly supervise students. In addition, you can split the necessary time needed to successfully bring a project to fruition, giving you more time for your other research.

My collaborator, Dr. Hoke, emphasized the importance of having structure for the students. With much success, we implemented the following plan for our summer students:

- Meet every morning with all students. Help them decide on the coming day's goals. Except when giving a lesson, try to keep these meetings short.
- For the first week, have materials ready for them to read and learn. If computer work is needed, have them begin learning how to use the necessary software. If you can incorporate the reading or computer work with beginning tasks of the project, all the better.
- Have the students present the material each week by giving a 10 minute talk. In the University of Richmond's Mathematics and Computer Science Department, these talks are given in a weekly meeting between all summer research groups. If you have multiple students, have two students present each week. Set up a schedule so they know who will present when and rotate the group members so that the students are not always working with the same person.
- If possible, meet again in the afternoon for a short period to see how work progressed that day and to review and begin setting goals for the next day.
- On Thursdays, have student presenters begin organizing the week's presentation into slides while the others continue to work on research.
- On Friday morning, review the slides and practice the talk.
- Have students write up their work as they go. Keeping track of their work in an organized and written way helps them solidify their knowledge and helps the project progress in a more methodical manner (and students also develop their writing skills). Publishing the work later is also easier when there is a written record of the project.

# 5  SELECT IMPACTFUL MENTORING EXPERIENCES

Here we discuss select mentoring experiences that had a significant impact on either our research program, or our approach to mentoring undergraduate research.

## 5.1 JRW's Experience at a Research Institution

Young faculty and postdoctoral fellows should make note of students who go above and beyond class assignments. These students are endogenously motivated and driven by more than just getting good grades in the classroom. These qualities make for good undergraduate researchers.

During my second semester as a postdoctoral fellow at Vanderbilt University, I taught a Mathematical Biology course. In that course, I instructed the students to explore current news sources and look for topics of interest to them that they might be able to mathematically model, and for them to derive questions that the models could address. Instead of reading the news, one of my students, Joe, read the scientific literature and melded ideas from current biotechnology with mathematical biology. His work went well beyond the assignment instructions, even including a rudimentary set of equations.

After the course ended, I asked Joe if he would like to continue working on the model as an undergraduate researcher under my supervision. Joe was an excellent student who could already read scientific literature and was a creative thinker. Working with him turned out to be a pleasure. After we developed the idea, I was able to bring in other postdoctoral researchers (one of them PSK) and we wrote a paper about the work that was published in the Journal of Biological Dynamics [7]. Afterward, we continued working on the project, finding an experimental collaborator in South Korea (Dr. Chae-Ok Yun) and then adding other scientists and mathematicians (including JLG). To date, we have four publications from these collaborations [7, 12, 23, 24].

Joe and I both learned a lot from working together on this project. I learned how to organize and develop a project for an undergraduate stepwise, breaking the work into weekly projects that he and I could accomplish. Joe was very busy with other work as well, and since the project was quite deep, I developed ideas about how to incorporate his work and ideas in the project while bringing other researchers in so that we could complete the more advanced parts of the project. Joe had a full research experience, bringing the project from conception to publication. I believe it was an essential experience that helped me land my current job at the University of Richmond, where faculty and administration are deeply interested in advancing undergraduate research. Our work also helped Joe get accepted into medical school, as beyond being a good student, Joe demonstrated that he was a creative scientific thinker.

**5.2 PSK's Experience at a Large Research University**

In 2012, within a year of arriving at the University of Sydney, I supervised James, who was an ideal first student, because he was independent and self-directed. For his project on cancer-immune modeling, I turned to two well-known papers [13, 14] and asked him to reproduce their results, extend the models, and conduct similar stability analyses of the extended system. I learned that a large part of supervising is finding the best way to motivate each individual, in this case challenging James immediately with advanced topics. James went on to pursue a PhD in applied mathematics at the University of Adelaide with Dr. Edward Green.

In 2013, I gave my next honors student Andrea a project stemming from one of my papers that presented an agent-based model (ABM) of cancer-immune dynamics. The goal was to formulate a partial differential equation (PDE) version of the ABM. One characteristic that stood out about Andrea was her perseverance to complete her thesis and finish with a high mark, despite circumstances that made her year much more difficult than normal. My most significant lesson from this project was that if students perceive that I care about them and believe in their potential, they push themselves beyond what I originally expected and find deep satisfaction in their accomplishments. My time with Andrea convinced me that I can give projects that are at the leading edge of my research and trust that students will produce something genuinely useful. Andrea's work was published the following year [6].

My increased faith in students' capacity made me more daring with projects. At the beginning of 2014, I met Jared and sent him away with two papers as starting points. He chose the more difficult one on modeling the evolution of human post-menopausal longevity. Like Andrea's project, the goal was also to study an ABM from one of my published papers and develop an analogous PDE system, but this time I did not know precisely where it would lead or what obstacles we might face. Jared seemed to catch the energy of the research, so although it was challenging, he dug into relevant papers, completed and even extended upon the original problem posed. Several people in the department said they thought it was the most interesting project that year. Jared received a scholarship to pursue a one-year master's degree with me, in the fall of 2015 he was admitted to the DPhil program in mathematics at Oxford University.

## 6 BEST PRACTICES AND PITFALLS: TWO SIDES OF THE SAME COIN

The easiest mistakes to make when mentoring undergraduate research students can be avoided by following best practices. One idea, introduced at the University of Richmond through a foundational program called The Long-Term Undergraduate Research Model (LURE) [1], helped the department become leaders in undergraduate mathematics research. LURE was funded by a grant from the National Science Foundation's Mentoring through Critical Transition Points (MCTP) program [1]. A core tenet of the program was to find students with research potential early in their college careers and supervise projects with them for multiple years. The idea is that students will learn the structure of the research process early and will have time in their later years to progress through a complete project. Through LURE, faculty members in the department learned how to develop a structure to summer research projects that best supports the needs of the nascent student researchers.

Beyond this, much has been written about best practices for undergraduate research [2, 5, 16, 19]. Here we will connect best practices with common pitfalls.

- *Choose the right project*: In interdisciplinary fields such as mathematical biology, one may feel compelled to give a student a pre-built model and just ask them to run simulations. While the student could make fairly quick progress on this front, this would not necessarily mean the student is learning much. After all, how is one learning to work like a mathematical biologist if the biology has already been studied, the model has already been formulated, and code has already been written by someone else? On the other hand, expecting the student to formulate a model from scratch, learn all the new mathematics to analyze the model, and learn how to write programs to simulate the model will often take more time than a student has. Choosing a problem where a preliminary model has been built and analyzed, but where a significant expansion is needed makes a great undergraduate project.

- *Expect Challenges*: Elsewhere, it has been written to ensure your students understand upfront that challenges are commonplace when doing research [5]. Beyond having your students expect such challenges, the faculty mentor needs to understand their role in helping the students deal with such challenges. On the one hand, you want the student to learn to persevere through these difficult moments, as this is what becoming a researcher will entail. On the other hand, you understand that the student only devotes so many hours to the project during the week, and that those hours can pass very quickly when a student is stuck alone with a problem (a particular challenge with single-student projects during the academic year). Having students document their issues and their attempts at solving the issue helps them to feel productive when stuck, but also allows them to know that the professor will be aware of their struggles soon enough. Learning to gently guide students away from a road block without fixing the problem for them is a skill that one acquires through practice. Be patient with yourself also.
- *Write it Down*: The importance of having your students write early and often cannot be understated, and elsewhere great advice for structuring a student's writing is given [5]. There, it is suggested that student work logs be split into four sections: hours worked each day, progress made, lingering questions, and additional thoughts/reflections. Even with this great structure, insisting on this without monitoring the student's writing is a recipe for this important step to become an afterthought for the student. Requiring students to email you their research log once a week (even if you don't actually have time to read it) will keep them on top of the process. More frequent updates may be required for full-time summer research. Plus, when you do have time to read it, it will help you to better understand your student's progress and struggles. In the end, writing often will help the student's prepare a final paper at the end of their project [2].

## 7 CONCLUSION

In this work, we have described our experiences working with undergraduate students on interdisciplinary research projects. For the projects to be successful, we find that students must have particular attributes, both in terms of their academic skill sets as well as their personalities. Additionally, support from the faculty member's institution is a factor in determining whether the faculty member and the research project will be successful. We also address best practices for forming and sustaining interdisciplinary undergraduate research groups. Finally, we discuss protocols for organizing individual as well as group projects.

Undergraduate research provides innumerable benefits to the students in terms of successful completion of college, admittance and success in graduate or professional school, career prospects, and growth as people in general [3, 8, 17, 18, 20]. For instance, a survey of 986 alumni from the University of Delaware found that those who had research experiences as undergraduates reported greater enhancement of cognitive skills and were more likely to pursue graduate degrees, among a host of other reported benefits [3]. Another study of undergraduate science research at nearly 100 institutions of higher education showed that such experiences force many students to take responsibility for their learning for the first time, encouraging levels of intellectual maturation that are often not achieved in a classroom setting [17].

Anecdotally, we have seen these benefits with our students manifest in many ways. As an example, undergraduate research has provided our students with the opportunity to write (or at least, co-write!) a scientific manuscript, and in some instances, learn about the peer review process and see their work published in journals and books [6, 7, 21]. The process of writing an article not only allows the student to reflect on their work, but also helps them place it in the context of pre-existing knowledge. Our students have also presented their research at local, national and even international conferences. This has given them excellent opportunities to enhance their presentation skills, while also putting them in a new environment to network with students and experts in their field.

Despite hurdles that may need to be overcome, undergraduate research can also result in rewards from a faculty perspective [22]. From our experiences, undergraduate research allows you to form deep relationships with your students as you watch them grow as thinkers and people. Beyond these personal/emotional benefits, we have also experienced scholarly benefits from mentoring undergraduate research projects. For instance, attending a conference (especially focused workshops) takes on a new perspective when you have undergraduate research students – we find we are often on a mission to find a vision for a current student's project and promising directions for future ones. The turnover of undergraduate research students keeps us returning to meetings with this exciting purpose in mind. As the years have progressed, undergraduate research has also allowed us to take on more projects than we could do without their support. Finally, undergraduate research can also open the door for new funding mechanisms with an undergraduate research mission [15], and can help you earn respect in your institution as your peers learn more about your research through your students. Lastly, undergraduate research projects can help you build collaborations within your institution, and even cross-institutionally, as was the case for us!